\documentclass[11pt,oneside,reqno]{amsart}

\usepackage{amssymb, amsmath, amsthm, mathtools, floatflt, palatino, euler, epic, eepic, floatflt, microtype}
\usepackage{cases}
\usepackage{graphicx}
\usepackage[usenames]{xcolor}

\usepackage{algorithm}
\usepackage[noend]{algpseudocode}

\makeatletter
\def\BState{\State\hskip-\ALG@thistlm}
\makeatother

\usepackage[margin=1 in]{geometry}
\definecolor{darkblue}{rgb}{0,0,0.7}
\definecolor{darkred}{rgb}{0.7,0,0}

\newcommand{\defi}[1]{\textit{\color{blue}#1}}  
\usepackage[colorlinks,linkcolor=darkred, citecolor=darkblue, pagebackref=true,pdftex]{hyperref}

\newcommand{\stan}{\operatorname{stan}}

\theoremstyle{plain}
\newtheorem{thm}{Theorem}[section]
\newtheorem{prop}[thm]{Proposition}
\newtheorem{lemma}[thm]{Lemma}
\newtheorem{cor}[thm]{Corollary}
\newtheorem*{theorem*}{Theorem}
\newtheorem*{cor*}{Corollary}
\newtheorem*{prop*}{Proposition}

\theoremstyle{definition}
\newtheorem{defn}[thm]{Definition}
\newtheorem{rem}[thm]{Remark}
\newtheorem{example}[thm]{Example}
\newtheorem{question}[thm]{Question}

\begin{document}

\title[Trees, parking functions, and standard monomials]{Trees, parking functions, and standard monomials \\of skeleton ideals}

\author{Anton Dochtermann}
\address{Department of Mathematics, Texas State University, San Marcos, TX, USA}
\email{dochtermann@txstate.edu}

 \author{Westin King}
\address{Department of Mathematics, Texas State University, San Marcos, TX, USA}
\email{westin\_king@alumni.baylor.edu}


\date{\today}


\begin{abstract}
Parking functions are a widely studied class of combinatorial objects, with connections to several branches of mathematics.   On the algebraic side, parking functions can be identified with the standard monomials of $M_n$, a certain monomial ideal in the polynomial ring $S = {\mathbb K}[x_1, \dots, x_n]$ where a set of generators are indexed by the nonempty subsets of $[n] = \{1,2,\dots,n\}$. Motivated by constructions from the theory of chip-firing on graphs we study generalizations of parking functions determined by $M^{(k)}_n$, a subideal of $M_n$ obtained by allowing only generators corresponding to subsets of $[n]$ of size at most $k+1$. For each $k$ the set of standard monomials of $M^{(k)}_n$, denoted $\stan_n^k$, contains the usual parking functions and has interesting combinatorial properties in its own right.

For general $k$ we show that elements of $\stan_n^k$ can be recovered as certain vector-parking functions, which in turn leads to a formula for their count via results of Yan. The symmetric group $S_n$ naturally acts on the set $\stan_n^k$ and we also obtain a formula for the number of orbits under this action. For the case of $k = n-2$ we study combinatorial interpretations of  $\stan_n^{n-2}$ and relate them to properties of uprooted trees in terms of root degree and surface inversions.  As a corollary we obtain a combinatorial identity for $n^n$ involving Catalan numbers, reminiscent of a result of Benjamin and Juhnke.  For the case of $k = 1$ we observe that the number of elements $\stan_n^1$ is given by the determinant of the reduced `signless' Laplacian, which provides a weighted count for $|\stan_n^1|$ in terms generalized spanning trees known as `spanning TU-subgraphs'.  Our constructions naturally generalize to arbitrary graphs and lead to a number of open questions.
\end{abstract}
\maketitle

\section{Introduction}\label{sec:intro}
A \defi{parking function of length $n$} is a sequence $a = (a_1, a_2, \dots, a_n)$ of nonnegative integers such that its rearrangement $c_1 \leq c_2 \leq \cdots \leq c_n$ satisfies $c_i < i$.  We let  ${\mathcal P}(n)$ denote the set of such sequences. This simple construction turns out to have connections and applications to many areas of mathematics including noncrossing partitions, hyperplane arrangements, and invariant theory (see \cite{YanPar} for a recent survey).  In \cite{KonWei} it is shown that the number of parking functions of length $n$ is given by $(n+1)^{n-1}$, which by Cayley's formula is equal to the number of spanning trees of the complete graph $K_{n+1}$. This correspondence can be generalized to the case of arbitrary graphs $G$ in the context of sequences of integers known as \emph{$G$-parking functions} \cite{Gab}.  For a general graph $G$ the number of $G$-parking functions is given by the number of spanning tree of $G$.  By the matrix tree theorem this number is given by $\det \tilde{\mathcal L}(G)$, the determinant of the reduced Laplacian of $G$.

On the algebraic side, the set ${\mathcal P}(n)$ can be identified with the standard monomials of the \defi{parking function ideal} $M_n$, a monomial ideal living in the polynomial ring $S = {\mathbb K}[x_1, x_2, \dots, x_n]$ for some fixed field ${\mathbb K}$.  This is the perspective adopted in \cite{PosSha}, where a parking function ideal $M_G$ is defined for an arbitrary graph $G$ on $n+1$ vertices.  Here we mostly restrict to the case where $G = K_{n+1}$ is the complete graph (see Section \ref{sec:conclusion} for a discussion of the general case). A generating set of monomials for the ideal $M_n$ is indexed by all nonempty subsets of $[n] = \{1,2,\dots, n\}$ according to the following construction.  For any nonempty $\sigma \subseteq [n]$  define the monomial $m_\sigma$ by 
\begin{equation}
m_\sigma =  \prod_{i \in \sigma} x_i^{n-|\sigma|+1}.
\end{equation}
The ideal $M_n$ is then by definition the ideal (minimally) generated by these monomials:
\[M_n = \langle m_\sigma: \emptyset \neq \sigma \subseteq [n] \rangle.\]
\noindent
The \defi{standard monomials} of $M_n$ are by definition the monomials which do \emph{not} appear in the ideal $M_n$ (which in turn form a basis for the ${\mathbb K}$-vector space $S/M_n$).  As discussed in \cite{PosSha} the standard monomials of $M_n$ correspond to the usual parking functions of length $n$.

Parking functions are closely related to the theory of \emph{chip-firing} on graphs, where an integer number of chips are placed on the vertices of a graph and passed to neighbors according to a simple `firing' rule.  The dynamics of this process are described by various notions of \emph{stability} of such configurations.  To recall some of these notions fix a root vertex and suppose ${\bf c} \in {\mathbb Z}^{n+1}$ is a configuration of chips with nonnegative values on the non-root vertices.  Then ${\bf c}$ is \defi{stable} if no individual non-root vertex can fire (the number of chips on any vertex is less than its degree) and is \defi{superstable} if no \emph{subset} of chips can fire simultaneously (without resulting in a vertex with a negative number of chips).  If $G$ is a graph on vertex set $\{0,1, \dots, n\}$ with specified root vertex $0$, one can check that the set of $G$-parking functions corresponds to the set of superstable configurations on the non-root vertices.  Such configurations are in a simple bijection with so-called \emph{critical} configurations, which form an abelian group called the \emph{critical group} of $G$.

In the context of chip-firing it is then natural to restrict the sets of vertices that are allowed to fire at once, interpolating between the notions of stable and superstable. Chip-firing with such restrictions has been studied in \cite{CarPaoSpo} and also \cite{Bac}, where the process was termed `hereditary chip-firing,' since the sets in question are required to be closed under taking subsets.   Motivated by these constructions, here we study subideals of $M_n$ generated by monomials corresponding to subsets of $[n]$ of a bounded size. Our main object of study will be certain \defi{$k$-skeleton ideals}, subideals of $M_n$ defined as follows.

\begin{defn}
For any $k = 0,1, \dots, n-1$, $M_n^{(k)}$ is the ideal in $S = {\mathbb K}[x_1, \dots, x_n]$ given by 
\[M_n^{(k)} = \langle m_\sigma: \sigma \subset [n],1 \leq |\sigma| \leq k+1 \rangle.\]
\end{defn}

Note that under this convention we have $M_n = M_n^{(n-1)}$. Also note that the generators of $M_n^{(k)}$ correspond to all subsets of size at most $k+1$, which can be thought of as the $k$-skeleton of a simplex.  Homological properties of these $k$-skeleton ideals are studied in \cite{DocOne}.  Here, we focus on combinatorial aspects.

In this context, the natural generalization of parking functions will be the standard monomials determined by the ideals $M_n^{(k)}$, by definition the set of exponent vectors of monomials that are not divisible by any element of $M_n^{(k)}$.  We let
\[\stan_n^k = \stan(M_n^{(k)})\]
denote the set of standard monomials of the ideal $M_n^{(k)}$.  Note that the usual parking functions are recovered as ${\mathcal P}(n) = \stan_n^{n-1}$.

\begin{example}
For $n = 4$ and $k=2$ we have
\[M_4^{(2)} = \langle x_1^4, x_2^4, x_3^4, x_4^4, x_1^3x_2^3, x_1^3x_3^3, x_1^3x_4^3, x_2^3x_3^3, x_2^3x_4^3, x_3^3x_4^3, x_1^2x_2^2x_3^2, x_1^2x_2^2x_4^2, x_1^2x_3^2x_4^2, x_2^2x_3^2x_4^2 \rangle.\]

The set $\stan_4^2$ consist of all monomials $x_1^{c_1}x_2^{c_2}x_3^{c_3}x_4^{c_4}$ which are not divisible by any of the generators.  In this case there are 152 such monomials (see Corollary \ref{cor:standard} below). We will typically think of these monomials as the sequences of integers given by the exponent vectors, so that the standard monomial $x_1^2x_2x_4^3$ is represented by $(2,1,0,3)$.
\end{example}

It turns out that set $\stan_n^k$ can be seen to coincide with certain \emph{vector parking functions}, introduced by Pitman and Stanley \cite{PitSta} and studied by Yan \cite{YanGen}.  Using results from \cite{YanGen} we obtain the following.


\newtheorem*{cor:standard}{Corollary \ref{cor:standard}}
\begin{cor:standard}
For any $n$ and $0 \leq k \leq n-1$, the number of standard monomials of $M_n^{(k)}$ is given by 
\[|\stan_n^k| = \sum \limits_{j=0}^k \binom{n}{j} (k + 1 - j)(k+1)^{j-1}(n-k)^{n-j}.\]
\end{cor:standard}

The set $\stan_n^k$ carries an action by the symmetric group $S_n$, and the orbits under this action can be seen to coincide with the set of (weakly) increasing elements of $\stan_n^k$.  In the case of parking functions, it is well-known that these orbits are counted by \emph{Catalan numbers}, and in the more general case of \emph{rational parking functions} they lead to a definition of the \emph{rational Catalan numbers} (see \cite{ArmLoeWar}).  In our case we obtain the following family of `skeletal Catalan numbers'.

\newtheorem*{thm:orbits}{Theorem \ref{thm:orbits}}
\begin{thm:orbits}
Let $S_n$ act on $\stan_n^k$ by permuting variables. Then the number of orbits under this action is given by
\[
\sum\limits_{\ell=1}^{k+1}\frac{\ell}{k+1}\binom{2k-\ell+1}{k}\binom{2(n-k)+\ell-2}{n-k+\ell-1}.
\]
\end{thm:orbits}

In the case of $k=n-1$ this sum leads to the decomposition of the Catalan numbers in terms of so-called \emph{ballot numbers}. We remark that in the case of $k=n-2$ it can also be shown that the number of orbits is given by the sum $C_n + C_{n-1}$ of consecutive Catalan numbers, where
\[
C_n = \dfrac{1}{n+1}\binom{2n}{n}.
\]
See Section \ref{sec:count} for details.

We next seek combinatorial interpretations of the set $\stan_n^k$, again in analogy with the study of parking functions.  Note that since $M_n^{(k)} \subset M_n^{(n-1)}$ for any $1 \leq k < n-1$ we have that $\stan_n^k$ strictly contains the set ${\mathcal P}(n)$ of parking functions.  Hence any combinatorial interpretation of $\stan_n^k$ should extend known interpretations of ${\mathcal P}(n)$.

We first study the case $k = n-2$.  We let $\partial (\stan_n^{n-2})$ denote the elements of $\stan_n^{n-2}$ that are not parking functions \footnote{In an earlier version of this paper posted on the arXiv the elements of $\partial(\stan_n^{n-2})$ were called `spherical parking functions'.}, so that
\[ \partial (\stan_n^{n-2}) := \stan_n^{n-2} \backslash \stan_n^{n-1}.\]
Note that an element of $\partial(\stan_n^{n-2})$ contains all variables $x_i$ with degree at least one. We can relate elements of $\stan_n^{n-2}$ to statistics on certain labeled trees.

\newtheorem*{thm:degree}{Theorem \ref{thm:degree}}
\begin{thm:degree}
For any $n \geq 1$ and $2 \leq s \leq n$, let $O(n,s)$ denote the number of elements in $\partial(\stan_{n}^{n-2})$ with $s$ variables of degree 1. Then $O(n,s)$ is given by the number of uprooted trees on $n$ having a root of degree $s-1$. This number is given by
\[O(n,s) = \binom{n}{s} (s-1)(n-1)^{n-s-1}.\]
It follows that $|\partial(\stan_n^{n-2})| = (n-1)^{n-1}$, the number of uprooted trees on $n$.
\end{thm:degree}

\begin{figure}[ht]
\includegraphics[scale = .475]{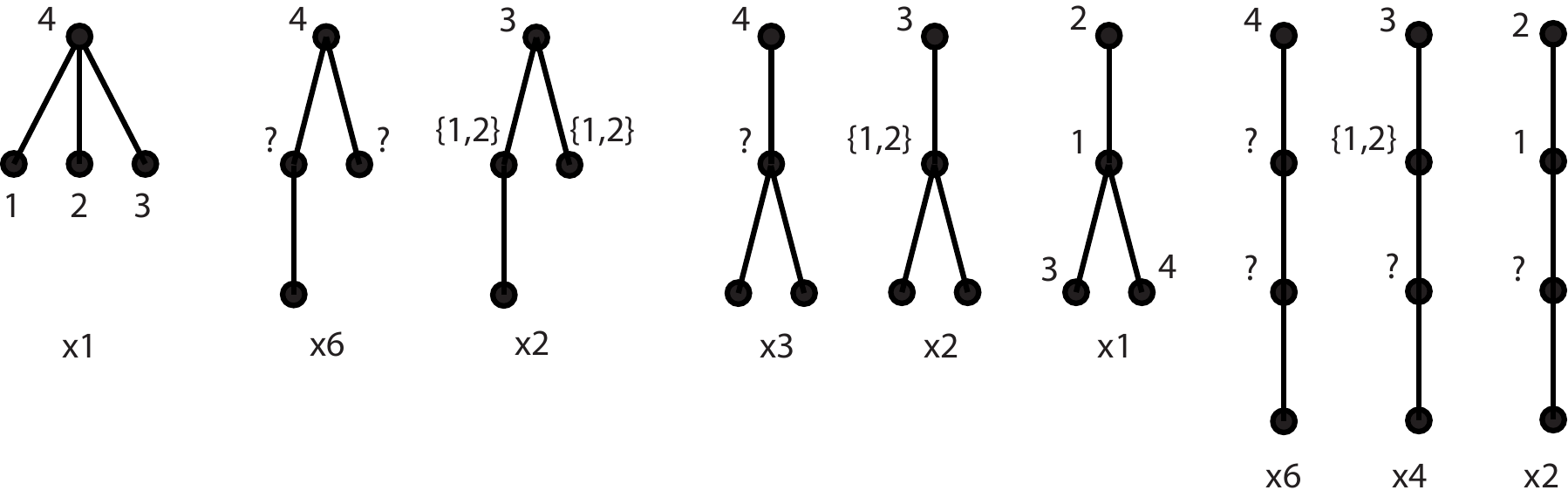}
\caption{Uprooted trees on $[4]$ by degree of the root, $27 = 1 + 8 + 18$.}
\label{fig:trees_degree}
\end{figure}

Here an \defi{uprooted tree on $n$} is a rooted tree with vertex set $[n]$ such that the root is larger than its `children' (immediate descendants). Figure \ref{fig:trees_degree} shows the uprooted trees on $[4]$.  The enumeration of uprooted trees (and other generalizations) is discussed in \cite{ChaDulGui}, and as an immediate corollary we get a simple count for the set $\stan_n^{n-2}$.

\newtheorem*{cor:stancodim}{Corollary \ref{cor:stancodim}}
\begin{cor:stancodim}
For any $n$ the number of standard monomials of $M_n^{(n-2)}$ is given by
\[|\stan_n^{n-2}| = (n+1)^{n-1} + (n-1)^{n-1}.\]
\end{cor:stancodim}

This of course follows from Corollary \ref{cor:standard} but we do not see a simple reduction.
After an earlier version of this paper was posted to the arXiv it was pointed out to the first author that the elements of $\partial(\stan_n^{n-2})$ can be seen to coincide with the set of \emph{prime parking functions} of length $n$ (after subtracting $1$ from each coordinate).  A \defi{prime parking function} (of length $n$) is a parking function of length $n$ that becomes a parking function of length $n-1$ when any 0 is deleted from the sequence. They were introduced by Gessel who also proved that the number of prime parking functions of length $n$ is equal to $(n-1)^{n-1}$ \cite{StaEC2}.  Our Corollary \ref{cor:stancodim} provides another proof of this fact.  This observation also indicates that prime parking functions are in fact the first step in a sequence of combinatorial objects coming from $\stan_n^{(k)}$. The construction also naturally generalizes to arbitrary graphs $G$, see Section \ref{sec:conclusion} for more discussion.

In the context of parking functions, a natural parameter to consider is \defi{degree}, by definition the sum of the entries. Kreweras \cite{Kre} showed that the number of elements of ${\mathcal P}(n)$ of degree $\binom{n}{2} - d$ is equal to the number of spanning trees of $K_{n+1}$ with $d$ \emph{inversions} (see Section \ref{sec:degree} for more details). A geometric interpretation of parking functions counted by degree is also provided by the so-called `Pak-Stanley' labeling of regions in the Shi arrangement \cite{Sta}.

In our context it turns out that the degrees of elements of $\partial(\stan_n^{n-2})$ are related to what we call \emph{surface inversions} (see Section \ref{sec:trees} for details).  We have the following result.

\newtheorem*{thm:inversions}{Theorem \ref{thm:inversions}}
\begin{thm:inversions}
Let $D(n,d)$ denote the number of elements of $\partial(\stan_n^{n-2})$ of degree $d$. Then $D(n,d)$ is equal to the number of uprooted trees on vertex set $[n]$ with $\binom{n}{2} - d +1$ surface inversions.
\end{thm:inversions}

The proofs of Theorem \ref{thm:degree} and Theorem \ref{thm:inversions} both involve modifications of Dhar's burning algorithm \cite{Dha} which provides a bijection between $G$-parking functions and spanning trees of an arbitrary graph $G$.  These are discussed in Section \ref{sec:burning}.

\begin{figure}[ht]
\includegraphics[scale = .475]{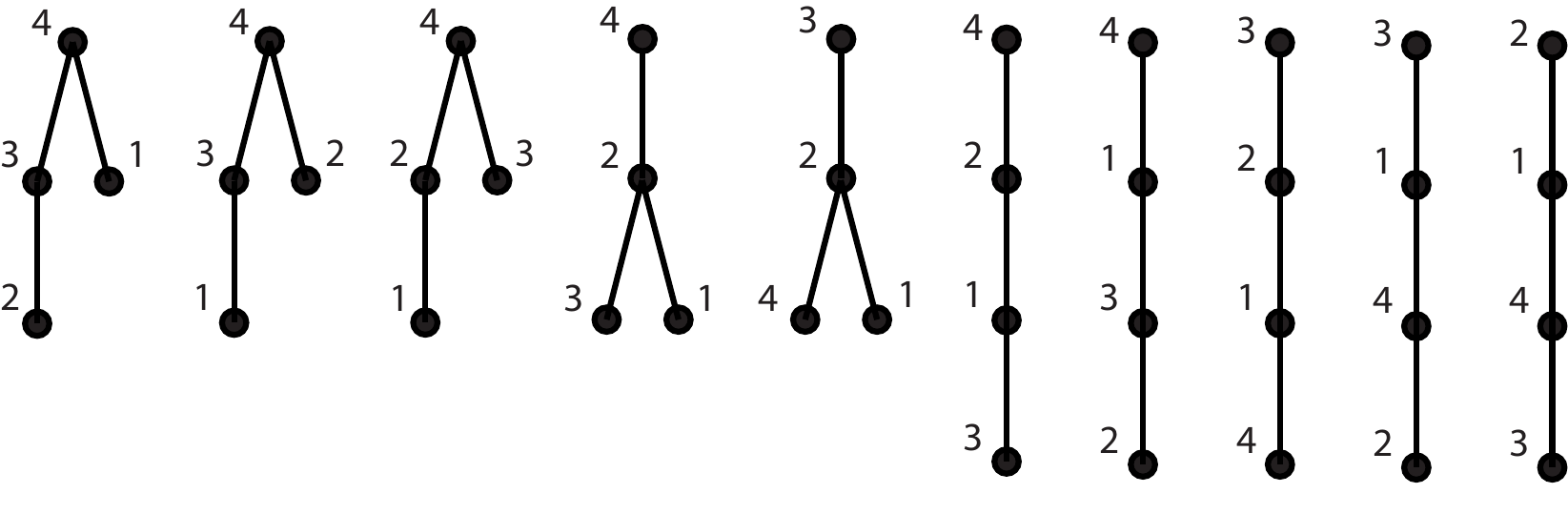}
\caption{Uprooted trees on $[4]$ with one surface inversion.}
\label{fig:trees_inversions}
\end{figure}

\begin{example}
We illustrate our results for the case of $n=4$ and $k=2$.  From Corollary \ref{cor:standard} we see that $\stan_4^2$ contains $5^3 + 3^3 = 152$ elements, and Theorem \ref{thm:orbits} says that these come in 19 $S_n$-orbits.  Among these, $125$ correspond to usual parking functions of length $4$ (for instance the parking function $(1,0,1,2)$ corresponds to the monomial $m = x_1x_3x_4^2$), which come in 14 orbits.  The other 27 standard monomials (the elements of $\partial(\stan_4^2))$ can be described by permuting variables of the following orbit representatives
\[x_1x_2x_3x_4, \hspace{.15  in} x_1x_2x_3x_4^2 \;(\times 4), \hspace{.15  in} x_1x_2x_3x_4^3\; (\times 4), \hspace{.15  in}x_1x_2x_3^2x_4^2\; (\times 6), \hspace{.15  in} x_1x_2x_3^2x_4^3\; (\times 12).\]
Note that if we count elements in terms of how many variables have degree $1$ we get $O(4,4) = 1$, $O(4,3) = 8$, and $O(4,2) = 18$. Theorem \ref{thm:degree} gives the identity
\[27 = 18 + 8 + 1.\]
Compare this to the degrees of the root in the uprooted trees depicted in Figure \ref{fig:trees_degree}.  

If, on the other hand, we count these elements by degree we see that $D(4,4) =1$, $D(4,5) = 4$, $D(4,6) = 10$, and $D(4,7) = 12$.  Theorem \ref{thm:inversions} then gives the identity
\[27 = 12 + 10 + 4 + 1.\]
See Figure \ref{fig:trees_inversions} for the $10$ uprooted trees with $\binom{4}{2} - 6 + 1 = 1$ surface inversion.
\end{example}

We next turn to the low dimensional skeleta.  First note that for $k=0$, we have that the ideal $M_n^{(0)} = \langle x_1^n, x_2^n, \dots, x_n^n\rangle$, so that we have $n^n$ standard monomials as predicted by Corollary \ref{cor:standard}. For $k=1$ the set $\stan_n^1$ turns out to have connections to other tree like structures via a determinantal interpretation. Recall that the usual parking functions (which correspond to the set $\stan_n^{n-1}$) are in bijection with the set of spanning trees of $K_{n+1}$, which by the Matrix Tree Theorem is counted by $\det(\tilde {\mathcal L}_{n+1})$, the determinant of the \emph{reduced Laplacian} of $K_{n+1}$.  For the case of $\stan_n^1$ a related matrix makes an appearance.

\newtheorem*{prop:signless}{Proposition \ref{prop:signless}}
\begin{prop:signless}
Let $\tilde {\mathcal Q}_{n+1} = \tilde {\mathcal Q}(K_{n+1})$ denote the reduced signless Laplacian of $K_{n+1}$. Then the number of standard monomials of $M_{n}^{(1)}$ is equal to the determinant $\det \tilde {\mathcal Q}_{n+1}$. This number is given by
\[|\stan_n^1| = \det \tilde {\mathcal Q}_{n+1} = (2n-1)(n-1)^{n-1}. \]
\end{prop:signless}

Via an application of the Cauchy-Binet theorem, one can see that $\det \tilde {\mathcal Q}_{n+1}$ has an interpretation as a weighted count of all \emph{spanning $TU$-subgraphs of $K_{n+1}$}, where a spanning $TU$-subgraph is a generalization of a spanning tree (bases of the matroid associated to a certain signed graph on $K_{n+1}$).  We refer to Section \ref{sec:signless} for details.   It would be interesting to find a bijection between $\stan_n^1$ and the set of all such graphs, ideally one that extends (any of the known) bijections between elements of $\mathcal{P}(n) = \stan_n^{n-1}$ and spanning trees.

The rest of the paper is organized as follows. In Section \ref{sec:vector} we discuss standard monomials of $M_G^{(k)}$ and relate them to vector parking functions. We also consider the orbits of the $S_n$-action on $\stan_n^k$.  In Section \ref{sec:trees} we focus on the case $k = n-2$ and show how elements of $\stan_n^{n-2}$ are related to the combinatorics of certain labeled trees via root degree and inversions.  In Section \ref{sec:signless} we consider the case $k = 1$ and relate the standard monomials of $M_n^{(1)}$ to the signless Laplacian and other tree-like structures. In Section \ref{sec:conclusion} we discuss our constructions in the context of general graphs and address some open questions.

\begin{rem}
While this version of the paper was being prepared, the preprint \cite{KumLatSon} was posted to the arXiv with some overlapping results.  In particular the authors make the connection to vector parking functions and give an independent proof of Theorem \ref{thm:inversions}.  They also describe how the ideals $M_n^{(k)}$ are related via Alexander duality to so-called \emph{multipermutohedron ideals}. 
\end{rem}

\section{Standard monomials, vector parking functions, and orbits}\label{sec:vector}

We briefly recall some basic notions in commutative algebra that we use throughout the paper.  We fix a field ${\mathbb K}$ and let $S = {\mathbb K}[x_1, \dots, x_n]$ denote the polynomial ring on $n$ variables with coefficients in ${\mathbb K}$.  If $M \subset S$ is any ideal generated by monomials, the \defi{standard monomials} of $M$ are by definition the set of monomials in $S$ that are not divisible by any element of $M$.   This notion makes sense for more general ideals $I \subset S$ where $S/I$ is finite dimensional as a ${\mathbb K}$-vector space (so-called Artinian ideals), in which case a ${\mathbb K}$-basis for $S/I$ is called a set of standard monomials.  

As mentioned above (and following the convention in \cite{PosSha}), any monomial $x_1^{a_1}x_2^{a_2} \cdots x_n^{a_n}$ can be identified with its exponent vector to provide a sequence $(a_1, a_2, \dots, a_n) \in {\mathbb N}^n$.  Hence when we discuss the set of standard monomials of an ideal $I$ we will often abuse notation and think about them as a collection of sequences, parking functions, etc.

Here we study the skeleton ideals $M_n^{(k)} \subset S$ defined in the introduction.  In particular $M_n^{(k)}$ is a monomial ideal generated by monomials $m_\sigma$, for all nonempty subsets $\sigma \subset [n]$ with $|\sigma| \leq k+1$.   We let $\stan_n^k$ denote the set of standard monomials of $M_n^{(k)}$.  By definition an element in $\stan_n^k$ is a monomial in $S$, but in analogy with parking functions we will often think of it as a sequence of nonnegative integers given by its exponent vector.  For example the element $x_1x_2 x_4^2x_5 \in M_5^{(2)}$ is represented by the sequence $(1,1,0,2,1)$.

It turns that these sequences are examples of \emph{vector parking functions}, a concept introduced by  Stanley and Pitman in \cite{PitSta} and further studied by Yan in \cite{YanGen}.  To recall the definition suppose ${\bf u} = u_1 \leq u_2 \leq \cdots \leq u_n$ is a sequence of weakly increasing positive integers.  A sequence $(a_1, a_2, \dots, a_n)$ of nonnegative integers is a \defi{${\bf u}$-parking function} if its rearrangement $c_1 \leq c_2 \leq \cdots \leq c_n$ satisfies $c_j < u_j$ for all $1 \leq j \leq n$. We will call  ${\bf u}$-parking functions and classical parking functions \defi{increasing} if $a_i \leq a_{i+1}$ for $1 \leq i \leq n-1$.

We let ${\mathcal P}({\bf u})$ denote the set of all ${\bf u}$-parking functions. Observe that the usual parking functions ${\mathcal P}(n)$ are recovered for the case ${\bf u} = (1,2, \dots, n)$.  We then have the following connection to our standard monomials.

\begin{lemma} \label{lem:vector}
For any $n \geq 1$ and $0 \leq k \leq n-1$ the set $\stan_n^k$ of standard monomials can be identified with the vector parking functions ${\mathcal P}({\bf u}_{n,k})$, where
\[{\bf u}_{n,k} = (\underbrace{n-k, n-k, \dots ,n-k}_\text{$n-k$ times}, n-k+1, n-k+2, \dots, n).\]
\end{lemma}

\begin{proof}
Recall that the ideal $M_n^{(k)}$ is generated by all monomials $m_\sigma$ where $\sigma \subseteq [n]$ and $1 \leq |\sigma | \leq k+1$.  Note that if $x_1^{a_1} x_2^{a_2} \cdots x_n^{a_n} \in \stan_n^{k}$, then no exponent can satisfy $a_i \geq n$ since $x_i^n$ is a generating monomial of $M_n^{(k)}$ for all $i$.  Hence $c_n < n$. 

If $k=0$ then $M_n^{(0)}$ is generated by the monomials $m_{\{i\}} = x_i^n$. Hence the only condition on the elements of $\stan_n^0$ is that every element must satisfy $c_i < n$.  We see that this corresponds to a ${\bf u}_{(n,0)}$-parking function.   Similarly, if $k=1$ then $M_n^{(1)}$ is generated by all monomials of the form $m_{\{i\}}$ and $m_{\{i,j\}} = x_i^{n-1}x_j^{n-1}$.  Hence an element of $\stan_n^1$ is characterized by having $c_n < n$ and at most one value $c_i$ satisfying $c_i = n-1$, in other words $c_i < n-1$ for all $i \leq n-1$.

Now suppose that the result holds for the elements of $\stan_n^{k-1}$, so that in particular $c_j < j$ for all $j \leq n-k+1$. Recall that $M_n^{(k-1)} \subset M_n^{(k)}$, and hence for general $k$ the elements of $\stan_n^k$ are characterized by the conditions of $\stan_n^{k-1}$ and the additional condition that at most $k$ values $c_i$ satisfy $c_i = n-k$, so that $c_i < n-k$ for all $i \leq n-k$.   The result follows.
\end{proof}

In \cite[Theorem 3]{Yan} Yan provides a formula for the number of ${\bf u}$-parking functions (for a certain class of vectors ${\bf u}$ that includes this case). Substituting $a=n-k$, $b=0$, $c=1$, and $m=k$ into her formula we obtain the following.

\begin{cor}\label{cor:standard}
The number of standard monomials of $M_n^{(k)}$ is given by 
\[|\stan_n^k| = \sum_{j=0}^k \binom{n}{j} (k + 1 - j)(k+1)^{j-1}(n-k)^{n-j}.\]
\end{cor}

Hence we have a good count on the number of elements in $\stan_n^{k}$.  In subsequent sections of the paper we seek combinatorial interpretations.

\subsection{Orbits}
For any $n \geq 1$ and $k = 0,1, \dots n-1$, the symmetric group $S_n$ acts on the set $\stan_n^k$ by permuting variables (this follows from the definition of $M_n$, which is evidently symmetric in the variables $x_i$).  One natural statistic to consider is the multiplicity of the trivial representation (the number of orbits) under this action. 

\begin{example}
Corollary \ref{cor:standard} tells us that the set $\stan_4^1$ consists of 189 elements.  The symmetric group $S_4$ acts on this set and we see that the number of orbits is given by $25 = 14 + 5 + 6$, corresponding to the \emph{increasing} representatives:
\begin{gather*}
\{0000,0001,0002,0003,0011,0012,0013,0022,0023,0111,0112,0113,0122,0123\} \\
\bigcup \; \{1111,1112,1113,1122,1123\} \\
\bigcup \; \{0222,0223,1222,1223,2222,2223\}.
\end{gather*}
\end{example}

For the case ${\mathcal P}(n) = \stan_n^{n-1}$ it is known that the number of orbits under the $S_n$-action is given by the Catalan number $C_n$, and in the context of \emph{rational parking functions} \cite{ArmLoeWar} the number of orbits is given by the so-called \emph{rational Catalan numbers} \cite{ArmLoeWar}  (see Remark \ref{rem:Dyck}).  Next we determine the number of orbits for the standard monomials of $\stan_n^k$ for all values of $k$.  

\begin{thm}\label{thm:orbits}
Let $S_n$ act on $\stan_n^k$ by permuting variables. Then, the number of orbits under this action is given by
\[
\sum\limits_{\ell=1}^{k+1}\frac{\ell}{k+1} \binom{2k-\ell+1}{k} \binom{2(n-k)+\ell-2}{n-k+\ell-1}.
\]
\end{thm}

\begin{proof}
We choose as our representative of each orbit the element with weakly-increasing exponent vector. From Lemma \ref{lem:vector} it is sufficient to count the increasing vector parking functions in the set $\mathcal{P}({\bf u}_{n,k})$. To accomplish this, we decompose any $c \in \mathcal{P}({\bf u}_{n,k})$ into two parts: an increasing classical parking function of length $k+1$ with $\ell$ instances of 0 and a weakly increasing sequence of length $n-k+\ell-1$ on the alphabet $\{0,1,\hdots,n-k-1\}$.

Let $c = (c_1,c_2,\hdots,c_n) \in \mathcal{P}({\bf u}_{n,k})$ satisfy $c_1\leq c_2\leq \hdots c_n$. We obtain a sequence that we can identify with a parking function of length $k+1$ by considering the vector $c' = (c_{n-k},c_{n-k+1},\cdots, c_n)$ then forming $d = (d_1,d_2, \hdots, d_{k+1})$ where $d_i = c_{n-k+i-1} - (n-k-1)$ when the difference is non-negative and 0 otherwise. Since $c_{n-k+i-1} < n-k + i-1$, we have $d_i < i$, so $d$ is indeed a classical increasing parking function. 

We are interested in the number of increasing classical parking functions of length $k+1$ with $\ell$ instances of 0, which can be enumerated via a straightforward application of the Cycle Lemma \cite{DvoMot} using the Dyck path representation of an increasing parking function, and is an entry in the Catalan Triangle (OEIS \#A009766):

\[
\dfrac{\ell}{k+1}\binom{2k-\ell+1}{k}.
\]

Now, we will count weakly-increasing sequences. Let $\ell$ be the largest index such that $d_i = 0$. By construction the numbers $c_1, \hdots, c_{n-k+\ell-1}$ form a weakly increasing sequence with the alphabet $\{0,1,\hdots,n-k-1\}$. The enumeration of such sequences is well-known (\cite{StaEC1}, section 1.2) and so we have 
\[
\binom{2(n-k)+\ell-2}{n-k+\ell-1}
\]
choices for our weakly increasing sequence.

Together, these two structures decompose our vector parking function $c$. On the other hand, given a parking function of length $k+1$ with $\ell$ instances of 0 and weakly increasing sequence of length $n-k+\ell-1$ on the alphabet $\{0,1,\hdots,n-k-1\}$, we can obtain an increasing vector parking function in $\mathcal{P}({\bf u}_{n,k})$ by adding $n-k-1$ to each entry in the parking function larger than 0 and appending those elements to our weakly increasing sequence.

Summing over all possible values for $\ell$ gives the result.
\end{proof}

\begin{example}
Consider the parking function $(0,0,0,2,4)$ and the sequence $(0,0,1,1,2)$. From these, we see $\ell = 3$, $k=4$, and thus $n = 5+4-3+1=7$. Since $n-k-1=2$, we add 2 to the non-zero entries in the parking function, then combine the two sequences to make the increasing vector parking function $(0,0,1,1,2,4,6) \in {\mathcal P}({\bf u}_{7,4})$.
\end{example}

\begin{rem} We may alternatively realize an increasing ${\bf u}$-parking function (those corresponding to our orbits) as a lattice path starting at $(0,0)$, ending at $(n,n-1)$ with steps $(1,0)$ and $(0,1)$ such that the path does not touch the points $\{(i-1,u_i)\}_{i=1}^n$. In this case, the number of such lattice paths (Equation 10.41, \cite{Kra}) is given by the determinant of a matrix whose entries are binomial coefficients:

\[
\det\limits_{1 \leq i,j \leq n} \left[ \binom{u_i}{j-i+1} \right].
\]
\end{rem}

\begin{rem}
Let $G_{n,k}$ denote the directed graph on vertex set $\{0,\hdots,n-1\}$ with edge set $\{(n-k-1,0)\} \cup \left( \bigcup\limits_{i=1}^{n-1} \{(i-1,i)\}\right)$. Then the elements in $\stan_n^k$ can also be realized as the parking functions of $G_{n,k}$, as defined in \cite{KingYan}. Furthermore, those monomials with weakly increasing exponent vectors can be identified with the increasing parking functions.
\end{rem}

\section{Uprooted trees, root degree, surface inversions, and $(n-1)^{n-1}$}\label{sec:trees}

In this section we specialize to the case $k=n-2$ and consider the combinatorics of $\stan_n^{n-2}$, the standard monomials of $M_n^{(n-2)}$. Once again we note that the ideal $M_n^{(n-2)}$ is contained in $M_n = M_n^{(n-1)}$ and hence we have $\stan_n^{n-1} \subset \stan_n^{n-2}$, the former set being the set ${\mathcal P}(n)$ of parking functions of length $n$ (of which there are $(n+1)^{n-1}$). Hence it is of interest to consider the `new contributions', for which we make the following definition.

\begin{defn}
For any integer $n \geq 1$ we let $\partial(\stan_n^{n-2}) = \stan_n^{n-2} \backslash \stan_n^{n-1}$.
\end{defn}
 
We can characterize the elements $\partial(\stan_n^{n-2})$ as follows.

\begin{lemma} \label{lem:rearrange}
 A sequence $a = (a_1, a_2, \dots, a_n)$ of nonnegative integers is the exponent vector of an element in $\partial (\stan_n^{n-2})$ if and only if its rearrangement $c_1 \leq c_2 \leq \cdots \leq c_n$ satisfies $c_1 = 1$  and $c_i < i$ for all $i = 2, 3,\dots, n$.
\end{lemma}

\begin{proof}
From Lemma \ref{lem:vector} we have that $\stan^{n-2}_n$ is given by the set of ${\bf u}$ parking functions, where ${\bf u} = (2,0,1,\dots, 1)$. Other the other hand the elements of $\stan^{n-1}_n$ are given by the ${\bf v}$-parking functions, where ${\bf v} = (1,1,\dots, 1)$. Note that both conditions imply that $c_i < i$  for all $i \geq 2$. Hence the elements of $\partial( \stan_n^{n-2})$ are characterized by the additional condition that $c_1=1$, and the result follows.
 \end{proof}

\begin{rem}\label{rem:prime}
A \defi{prime parking function} of length $n$ is a parking function of length $n$ that becomes a parking function of length $n-1$ after removing any $0$ entry \cite[Exercise 5.49f]{StaEC2}. From Lemma \ref{lem:rearrange} we see that $\partial(\stan_n^{n-2})$ can be seen to coincide with the set of prime parking functions (after subtracting 1 from each coordinate).  Hence the set $\stan_n^{n-2}$ contains two disjoint sets $A$ and $B$, where $A$ is the set of monomials corresponding to parking functions and $B$ is in (easy) bijection with the set of prime parking functions of length $n$.
\end{rem}

\begin{rem}\label{rem:Dyck}
The set $\partial(\stan_n^{n-2})$ can also be seen to be a special case of \emph{rational parking functions} as discussed in \cite{ArmLoeWar}.   In their language the elements of $\partial(\stan_n^{n-2})$ are $(n, n-1)$-parking functions, and correspond to lattice paths weakly staying above a line of slope $n/(n-1)$ (so-called \emph{rational Dyck paths}).  By results in \cite{ArmLoeWar} the set of such sequences has cardinality $(n-1)^{n-1}$, and carries a permutation representation of $S_n$ where the number of orbits is given by the (usual) Catalan number $C_{n-1}$.  The decomposition of $(n-1)^{n-1}$ given in Corollary \ref{cor:another} describes this set in terms of $S_n$-orbits. 
\end{rem}

\subsection{Breadth-first and depth-first burning algorithms} \label{sec:burning}
We wish to relate the elements of $\stan_n^{n-2}$ to certain statistics on labeled trees. For this we will employ `depth-first search' (DFS) and `breadth-first search' (BFS) variants of the classical \emph{burning algorithm} due to Dhar \cite{Dha}.

The `depth-first' version of the algorithm is studied by  Perkinson, Yang, and Yu \cite{PerYanYu} and may be informally stated as follows. Suppose $p = (p_1, \dots, p_n)$ is any sequence of nonnegative integers. Consider the graph $K_{n+1}$ with root 0 and $p_i$ number of firefighters standing on vertex $i>0$. A fire is lit at 0 and attempts to spread to the vertex $n$. If $p_n > 0$, then one of the firefighters wets the edge $\{0,n\}$ and the fire fails to spread. The fire then attempts to spread to $n-1$, following the same procedure. The fire attempts to spread from its current position to the largest adjacent un-burnt vertex and $p_i$ edges are wet by firefighters before the fire can spread to vertex $i$. If the fire can not spread from its current vertex, it retreats to the parent vertex and attempts to spread to the largest adjacent vertex along a non-wet edge.  In \cite{PerYanYu} it is shown that $p$ is a parking function if and only if the collection of burnt edges forms a spanning tree of $K_{n+1}$.  The following is Algorithm 1 from \cite{PerYanYu}, modified for our purposes, and applied to the graph $K_{n+1}$.

\begin{algorithm}
\caption{DFS-burning algorithm}\label{alg:DFS}
\begin{algorithmic}[1]
\Statex \textbf{Input}: $p: \in \mathbb{N}^n$
\BState \texttt{burnt\_vertices} $=\{0\}$
\BState \texttt{tree\_edges} $=\{\}$
\BState execute \textsc{dfs\_from}$(0)$
\Statex \textbf{Output}: \texttt{tree\_edges}
\Statex
\Function{dfs\_from}{$i$}
    \ForAll{$j \in [n]$, from largest to smallest}
        \If{$p_j =0$}
            \State append $j$ to \texttt{burnt\_vertices}
            \State append $\{i,j\}$ to \texttt{tree\_edges}
            \State $p_j = p_j-1$
            \State \textsc{dfs\_from}$(j)$
        \Else
            \State $p_j = p_j-1$
        \EndIf
    \EndFor
\EndFunction
\end{algorithmic}
\end{algorithm}

We will also need a `breadth-first' version of the algorithm, which proceeds similarly to the above except for the fact that the fire attempts to burn `level by level'.  The details are given in Algorithm \ref{alg:BFS}.  Since we could not find this particular process described in the literature, we prove here that in fact it gives the desired bijection.

\begin{lemma}
Algorithm \ref{alg:BFS} provides a bijection between the set ${\mathcal P}(n)$ of parking functions of length $n$ and the set of spanning trees of $K_{n+1}$.
\end{lemma}

\begin{proof}
We first argue that Algorithm \ref{alg:BFS} returns a spanning tree if and only if the input was a parking function. For each vertex added to \texttt{burnt\_vertices} after initiation, one edge incident to that vertex and one other in \texttt{burnt\_vertices} is added to \texttt{tree\_edges}. The resulting graph is connected, so if all vertices are in \texttt{burnt\_vertices} at termination, then the result must be a spanning tree.

Suppose for $p \in \mathbb{N}^n$ that the algorithm does not terminate in a spanning tree. Then for each $i \notin \texttt{burnt\_vertices}$, we must have $p_i \geq |\texttt{burnt\_vertices}|$. However, this means $p$ is not a parking function since after rearrangement, $c_{|\texttt{burnt\_vertices}|+1} \geq |\texttt{burnt\_vertices}|+1.$

On the other hand suppose for $c \in \mathbb{N}^n$ that the algorithm does terminate in a spanning tree and also assume without loss of generality that $c_i \leq c_{i+1}$ for $1 \leq i < n$. Notice that the vertices will be added to \texttt{burnt\_vertices} in ascending order, meaning $c_i - j = 0$ for some $j < i$, which in turn means $c_i < i$, making $c$ a parking function.

We now argue that Algorithm \ref{alg:BFS} is injective. Given two distinct length $n$ parking functions, $s$ and $p$, and for an appropriate $i$, assume without loss of generality that $s_i$ is minimal among all $s_j$ and $p_j$ such that $s_j \neq p_j$. Then the algorithm will add the same vertices to \texttt{burnt\_vertices} for both $s$ and $p$ until the point that $s_i$ is reduced to zero (and $p_i$ is not). On the next iteration of the `for' loop in line 7, an edge incident to vertex $i$ will be added in the case of $s$, but the same edge will not be added in the case of $p$, causing the resulting trees to be distinct. Finally, the algorithm must give a bijection because the number of spanning trees of $K_{n+1}$ and the number of length $n$ parking functions are the same.

We can easily reverse Algorithm \ref{alg:BFS} as follows: let $T$ be a tree spanning $K_{n+1}$ and for $0 \leq i \leq n$, let $i \in L_j$ if and only if there are $j$ edges on the shortest path between 0 and $i,$ and let $\text{par}(i)$ be the parent of $i$ in $T$. Then for $i>0$, if $i \in L_j$, we have $p_i = |\{m \in L_{j-1} : m < \text{par}(i)\}|+ \sum_{\ell=0}^{j-2}|L_\ell|$.
\end{proof}

\begin{algorithm}
\caption{BFS-burning algorithm}\label{alg:BFS}
\begin{algorithmic}[1]
\Statex \textbf{Input}: $p \in \mathbb{N}^n$
\BState \texttt{burnt\_vertices} $=\{0\}$
\BState \texttt{tree\_edges} $=\{\}$
\BState \texttt{current\_level} $= \{0\}$
\BState \texttt{next\_level} $=\{\}$
\While{\texttt{current\_level}$\neq \{\}$}
    \ForAll {$i \in$ \texttt{current\_level}, from smallest to largest}
        \ForAll{$j \in [n]$}
            \If{$p_j =0$}
                \State append $j$ to \texttt{burnt\_vertices}
                \State append $j$ to \texttt{next\_level}
                \State append $\{i,j\}$ to \texttt{tree\_edges}
                \State $p_j = p_j-1$
            \Else
                \State $p_j = p_j-1$
            \EndIf
        \EndFor
    \EndFor
    \State \texttt{current\_level} $=$ \texttt{next\_level}
    \State \texttt{next\_level} $=\{\}$
\EndWhile
\Statex \textbf{Output}: \texttt{tree\_edges}
\end{algorithmic}
\end{algorithm}

\subsection{Root degrees and surface inversions}

Next we recall some results from the theory of tree enumeration. Let ${\mathcal T}_{n,k}$ denote the family of rooted labeled trees on $[n]$ such that the root has exactly $k$ larger immediate descendants.  We will refer to elements in ${\mathcal T}_{n,0}$ as \defi{uprooted trees}.  In \cite{ChaDulGui} the authors provide formulas for the size of the sets ${\mathcal T}_{n,k}$ counted by various statistics, and in particular establish the following.

\begin{prop}[\cite{ChaDulGui}, Proposition 4.2] \label{prop:treesdegree}
The number of uprooted trees on $[n]$ with root degree $s-1$ is given by 
\[\binom{n}{s} (s-1) (n-1)^{n-s-1}.\]
From this it follows that 
\[(n-1)^{n-1} = \sum_{s=1}^{n-1} \binom{n}{s} (s-1) (n-1)^{n-s-1}.\]
\end{prop}

We next use these observations to study elements of $\partial(\stan_n^{n-2})$.

\begin{thm} \label{thm:degree}
For any $n \geq 1$ and $2 \leq s \leq n$  let $O(n,s)$ denote the number of elements in $\partial(\stan_{n}^{n-2})$ with $s$ variables of degree 1. Then $O(n,s)$ is given by the number of uprooted trees on $[n]$ having a root of degree $s-1$. This number is given by
\[O(n,s) = \binom{n}{s} (s-1)(n-1)^{n-s-1}.\]
The total number of elements in $\partial(\stan_n^{n-2})$ is given by
\[|\partial(\stan_n^{n-2})| = (n-1)^{n-1}.\]
\end{thm}

\begin{proof}
We utilize the BFS-burning algorithm, Algorithm \ref{alg:BFS}, to provide a bijection between the elements of $\partial(\stan_n^{n-2})$ with $s$ variables of degree 1 and uprooted trees on the vertex set $[n]$ in which the root has $s-1$ children. To this end, let $a$ be such an element of  $\partial(\stan_n^{n-2})$. Per Remark \ref{rem:prime}, we obtain a prime parking function $p$ by letting $p_i = a_i-1$ for $i \in [n]$. Furthermore, we also know we may delete any 0 from $p$ to obtain a parking function of length $n-1$, so let us delete the final $p_j=0$ in $p$ to obtain $\hat{p}$ and we remember $j$ for later. We note that there are $(s-1)$ instances of 0 in $\hat{p}$.

Now, we run Algorithm $\ref{alg:BFS}$ on $\hat{p}$ and we notice that the root will have child $i$ if and only if $\hat{p}_i = 0$. Thus, our resulting tree's root has $s-1$ children, but is not an uprooted tree as the root is labeled 0. This is easily fixed by relabeling the root by $j,$ the index of the previously-deleted 0, and increasing all non-root vertex labels by 1 if they are $j$ or larger. Our tree's root is now larger than its children. That this process is indeed a bijection follows easily from the fact that Algorithm \ref{alg:BFS} is reversible.

Hence, by Proposition \ref{prop:treesdegree}, we have found that $O(n,s) = \binom{n}{s} (s-1)(n-1)^{n-s-1}$ as claimed.
\end{proof}

\begin{rem}
One can see that the argument used in the proof of Theorem \ref{thm:degree} provides an analogous statement for the case of ordinary parking functions.  That is, the number of elements in ${\mathcal P}(n)$ with $s$ variables equal to 0 is given by the number of rooted trees on $\{0,1,...,n\}$ with root 0 where the root has degree $s-1$.  This also follows from a result in \cite{ChePyl}, where the authors describe a bijection due to Postnikov from the Dyck path representation of the elements of ${\mathcal P}(n)$ to trees rooted at 0. 
\end{rem}

From Theorem \ref{thm:degree} we get a formula for the number of elements in $\stan_n^{(n-2)}$.  The following can presumably be derived from Corollary \ref{cor:standard} but we do not see a simple reduction.  

\begin{cor}\label{cor:stancodim}
The number of standard monomials of $\stan_n^{(n-2)}$ is given by
\[|\stan_n^{n-2}| = (n+1)^{n-1} + (n-1)^{n-1}.\]
In particular we have
\[|\partial(\stan_n^{n-2})| = (n-1)^{n-1}.\]
\end{cor}

\begin{proof}
By definition we have that $\stan_n^{n-2} = \stan_n^{n-1} \cup \; \partial(\stan_n^{n-2})$ (disjoint union). Recall that $\stan_n^{n-1} = {\mathcal P}(n)$ has cardinality $(n+1)^{n-1}$.  The result follows from Theorem \ref{thm:degree}.
\end{proof}

\begin{rem}
From Remark \ref{rem:prime} we saw that the elements of $\partial(\stan_n^{n-2})$ can be identified (after subtracting one from each entry) with the set of prime parking functions of length $n$. Hence Theorem \ref{thm:degree} provides another proof that there are $(n-1)^{n-1}$ such sequences, a fact that was first established by Gessel \cite{StaEC2}.
\end{rem}

\subsection{Degree and surface inversions}\label{sec:degree}

We next turn to a proof of Theorem \ref{thm:inversions}.  Recall that in the classical case Krewaras \cite{Kre} showed that the number of length $n$ parking functions of degree ${n \choose 2} - k$ is equal to the number of spanning trees of $K_{n+1}$ with $k$ \emph{inversions}.  To describe this notion suppose $T$ is a tree on vertex set $\{1, 2,\dots,n\}$ and fix a root at 1.  Then an \defi{inversion} in $T$ is a pair $(\alpha, \beta)$ of vertices in $T$ such that $\beta$ is a descendant of $\alpha$ and $\alpha > \beta$. Kreweras established his result by relating parking functions to the \emph{Tutte polynomial} of $K_n$, and hence to the \emph{external activity} of trees.  A bijection between the set of trees with $k$ inversions and the set of trees with $k$ externally active edges was given by Beissinger in \cite{Bei}.  A geometric interpretation of parking functions counted by degree is also provided by the so-called `Pak-Stanley' labeling of regions in the Shi arrangement \cite{Sta}.  Under this labeling, the set of prime parking function correspond to the \emph{bounded} regions.

In our context we can equally well enumerate elements of $\partial(\stan_n^{n-2})$ according to degree.  As we have seen, elements in $\partial(\stan_n^{n-2})$ are in bijection with the set of \emph{uprooted trees} on $n$, and so degree gives a new parameter on such objects.   For this we need the following notion.

\begin{defn}
Suppose $T$ is a uprooted tree on $[n]$ with root $r$. A \defi{surface inversion} is a pair of \emph{non-root} vertices $(\alpha, \beta)$ such that $\beta$ is a descendant of $\alpha$ and $\alpha > \beta$. 
\end{defn}

We refer to Figure \ref{fig:trees_inversions} for examples. We then have the following correspondence.

\begin{thm}\label{thm:inversions}
Let $D(n,d)$ denote the number of elements in $\partial(\stan_n^{n-2})$ of degree $d$. Then $D(n,d)$ is equal to the number of uprooted trees on vertex set $[n]$ with ${n \choose 2} - d + 1$ surface inversions.
\end{thm}

\begin{proof}
We will employ Algorithm \ref{alg:DFS} to prove the claim via bijection. Theorem 3 of \cite{PerYanYu} concludes that Algorithm \ref{alg:DFS} gives a bijection between ${\mathcal P}(n)$ and spanning trees of $K_{n+1}$ such that the image of a parking function $p$ has ${n \choose 2} - \deg(p)$ inversions.

Now consider $s \in \partial(\stan_n^{n-2})$ of degree $\deg(s) = d$. Per Remark \ref{rem:prime}, we may subtract 1 from each entry to obtain the prime parking function $p = (p_1,p_2,\hdots,p_n)$. Recall that if any 0 is removed from a prime parking function of length $n$, then the resulting sequence is a parking function of size $n-1$. Delete the final 0, $p_i$ for some $i$, from $p$ and consider the resulting sequence $\hat{p} = (\hat{p}_1,\hat{p}_2,\hdots,\hat{p}_{n-1})$ where $\hat{p}_j = p_j$ for $j < i$ and $\hat{p}_j = p_{j+1}$ for $j \geq i$. We note that $\deg(\hat{p}) = deg(p)$. Now run the Algorithm \ref{alg:DFS} on the complete graph $K_n$. The resulting spanning tree has ${n-1 \choose 2} - \deg(\hat{p}) = {n-1 \choose 2} - (\deg(s) - n) = {n \choose 2} - d + 1$ inversions. Since $\hat{p}_j > 0$ for $j \geq i$, all of the root's children must be smaller than $i$ as the edges $\{0,j\}$ were wet before the fire first spread from the root during the burning algorithm. Furthermore none of the inversions involve the root as the root has label 0.

We finish by converting $T$ into an uprooted tree by relabeling the vertices $j = j+1$ for $j \geq i$ then relabeling the root with $i$. The ${n \choose 2} - d + 1$ inversions before relabeling are now surface inversions as they do not involve the root and all of the root's children have labels less than $i$ since the root before relabeling had no children with labels $i$ or larger.

One can see that the above process is reversible. Assuming we start from an uprooted tree with vertex set $[n]$ and ${n \choose 2} - d + 1$ surface inversions we convert it to a spanning tree of $K_n$ by relabeling the root 0 and adjusting the other labels appropriately. Apply Algorithm \ref{alg:DFS} in reverse to get a parking function of length $n-1$.  This is then converted into our desired element of $\partial(\stan_n^{n-2})$ by inserting a 0 at the coordinate given by the original root's label and then adding 1 to all entries.
\end{proof}

\begin{figure}[ht]
\includegraphics[scale = .55]{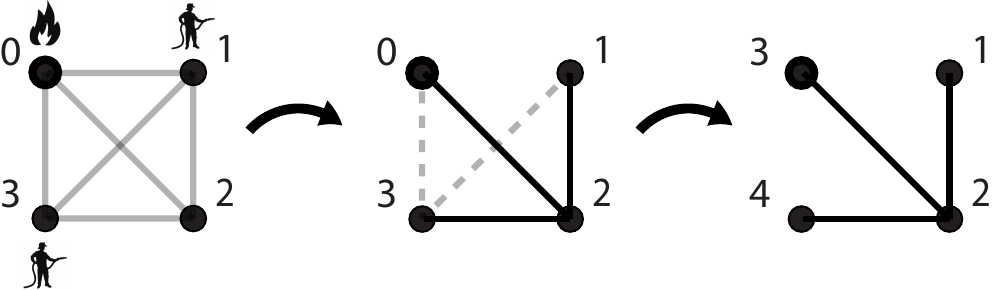}
\caption{Obtaining an uprooted tree from $s = (2,1,1,2) \in \partial(\stan_4^2)$.}
\label{fig:burning}
\end{figure}

\begin{example}
We give an example of the bijection described above with $n=4$. Let $s=(2,1,1,2),$ so that $p=(1,0,0,1)$ and $\hat{p}=(1,0,1)$. We have $\deg(s)=6$, and so we expect an uprooted tree with ${4 \choose 2} - 6 + 1=1$ surface inversions. The final 0 of $p$ is in the 3rd position, so $i=3$. 

Running Algorithm \ref{alg:DFS} with $\hat{p}$ on $K_4$ with vertices $\{0,1,2,3\}$, the fire attempts to spread to vertex 3 but the edge $\{0,3\}$ is dampened. The fire spreads to vertex 2, then to vertex 3 (as one edge was previously dampened). The fire next attempts to spread to vertex 1 but the edge $\{1,3\}$ is dampened and the fire finally spreads to vertex 1 via vertex 2.

Since $i=3$, we relabel all vertices with labels $3$ or larger and then relabel the root as 3. The resulting tree has $1$ surface inversion and the root is larger than its child. Figure \ref{fig:burning} shows $K_4$ before the algorithm initializes, after the burning algorithm runs (where dashed edges are those dampened), as well as the tree after relabeling.
\end{example}

\begin{rem}
We have seen that the set $\stan_n^{n-2}$ forms a ${\mathbb K}$-basis for the $S$-algebra $S/M_n^{(n-2)}$, and hence the degree sequence of $\stan_n^{n-2}$ provide the coefficients for the Hilbert series of the finite dimensional ${\mathbb K}$-vector space $S/M_n^{(n-2)}$.  We do not know if the subset $\partial(\stan_n^{n-2})$ can be seen as the set of standard monomials for some naturally occurring $S$-algebra.  If this was the case, Theorem \ref{thm:inversions} would provide a combinatorial interpretation for the coefficients of its Hilbert series.
\end{rem}

\subsection{Another way to count $(n-1)^{n-1}$}\label{sec:count}

We next describe how an explicit formula for $|\stan^{n-2}_n|$ coming directly from the definition provides another way of counting $(n-1)^{n-1}$ in the spirit of \cite{BenJuh}. 

\begin{cor} \label{cor:another}
For any integer $n \geq 1$ we have 
\begin{equation} \label{eq:choose}
(n-1)^{n-1} \; = \sum\limits_{\substack{0 \leq k_1 \leq 1 \\ 0 \leq k_1 + k_2 \leq  2 \\ \cdots \\ 0 \leq k_1 + \cdots + k_{n-2} \leq n-2}} {n \choose k_1}{n-k_1 \choose k_2} \cdots {n-(k_1 + \dots + k_{n-3}) \choose k_{n-2}},
\end{equation}
where $n > 1$ is an integer and $k_1, k_2, \dots, k_{n-2}$ are nonnegative integers.
\end{cor}

\begin{proof}
Our strategy will be to show that the right hand side of the identity naturally counts the elements of $\partial(\stan_n^{n-2})$, Corollary \ref{cor:stancodim} will then give the result. Note that if $x_1^{a_1} x_2^{a_2} \cdots x_n^{a_n} \in \partial(\stan_n^{n-2})$, then no exponent can satisfy $a_i \geq n$ since $x_i^n$ is a generating monomial of $M_G^{(n-1)}$.  We can have at most one exponent satisfying $a_i = n-1$, at most two exponents satisfying $a_j = n-2$, and so on.  Let $k_1$ denote the number of exponents $a_i$ such that $a_i = n-1$, $k_2$ the number of exponents $a_j$ such that $a_j = n-2$, etc.  Hence the number of standard monomials of the form $x_1^{a_1} x_2^{a_2} \cdots x_n^{a_n}$ is given by the expression on the right hand side of (\ref{eq:choose}), which by Corollary  \ref{cor:stancodim} is equal to $(n-1)^{n-1}$.
\end{proof}

Note that the expression in the summand in Corollary \ref{cor:another} can be written as $\binom{n}{k_1,k_2, \dots, k_{n-2},k_{n-1}}$, where $k_{n-1} = n - (k_1 + \cdots + k_{n-2})$.  This allows for more direct comparison with the identity in Equation \ref{eq:another} below. To illustrate the identity, for $n=4$ we have the possible values of $k_1$, $k_2$ given by
\begin{center}
\begin{tabular}{|c|c|}
\hline
$k_1$ & $k_2$ \\
\hline
0 & 0 \\
0&1 \\
0&2 \\
1&0 \\
1&1 \\
\hline
\end{tabular}
\end{center}

and Equation \ref{eq:choose} becomes
\smallskip
\begin{align*}
3^3 &= {4 \choose 0}{4-0 \choose 0} + {4 \choose 0}{4-0 \choose 1} + {4 \choose 0}{4-0 \choose 2} + {4 \choose 1}{4-1 \choose 0} + {4 \choose 1}{4-1 \choose 1} \\[8 pt]
&= 1 + 4 + 6 + 4 + 12.
\end{align*}

Note that the formula for $(n-1)^{n-1}$ involves a total of $C_{n-1}$ terms in the summation.  In addition, each summand represents an orbit in the $S_n$-action on the elements of $\partial(\stan_n^{n-2})$.  Hence for the case $k = n-2$ we see that the formula from Theorem \ref{thm:orbits} for the number of orbits of $\stan_n^{n-2}$ reduces to $C_n + C_{n-1}$.

In \cite{BenJuh} Benjamin and Juhnke established a similar looking identity:

\begin{equation} \label{eq:another}
(n-1)^{n-1} \; = \sum\limits_{\substack{0 \leq k_1 \leq 1 \\ 0 \leq k_1 + k_2 \leq  2 \\ \cdots \\ 0 \leq k_1 + \cdots + k_{n-2} \leq n-2}} \frac{(n-1)!}{k_1!k_2! \cdots k_{n-2}!},
\end{equation}

\noindent
where $n > 1$ is an integer and $k_1, k_2, \dots, k_{n-2}$ are nonnegative integers.  In \cite{DWC} this formula was generalized to an identity involving $n^m$, where $m<n$.  Note that Equation \ref{eq:another} involves a summation over the same indexing set as Equation \ref{eq:choose} (of size given by a Catalan number), but one can check that the terms in the summation are not the same (even up to reordering).  For instance if $n=4$ Equation \ref{eq:another} becomes

\begin{align*}
3^3 &= \frac{3!}{0!0!}+\frac{3!}{0!1!}+ \frac{3!}{0!2!}+ \frac{3!}{1!0!} + \frac{3!}{1!1!} \\[8 pt]
&= 6 + 6 + 3 + 6 + 6.
\end{align*} 

We do not know if there is a `parking function' interpretation of Equation \ref{eq:another}.

\section{Signless Laplacians and TU-subgraphs}\label{sec:signless}

In this section we discuss combinatorial interpretations of $\stan_n^1$, the standard monomials of the ideal $M_n^{(1)}$.  Recall that $M_n^{(1)}$ is generated by all $m_\sigma$ with $\sigma \subset [n]$ and $1 \leq |\sigma| \leq 2$. As we have seen, for any $n \geq 1$, the ideal $M_n^{(1)}$ is a subideal of $M_n$ and hence the set $\stan_n^1$ contains the set $\stan_n^{n-1}$ (which correspond to parking functions of length $n$).

\begin{example}\label{ex:K4oneskel}
For $n=3$ we have $M_3^{(1)} = \langle x_1^3, x_2^3, x_3^3, x_1^2x_2^2, x_1^2x_3^2, x_2^2x_3^2 \rangle$, with the set of standard monomials given by 
\[\stan_3^2 \cup \{x_1x_2x_3, x_1^2x_2x_3, x_1x_2^2x_3, x_1x_2x_3^2\},\]
\noindent
giving a total of $16+4 = 20$ standard monomials.   Here as usual we identify a parking function with the monomial having that sequence as an exponent vector, so that for example $(1,0,2)$ is identified with the monomial $x_1x_3^2$.
\end{example}

From Proposition \ref{cor:standard}  we get a formula for the number of elements in $\stan_n^1$.   We also provide a direct proof.

\begin{prop}\label{prop:standardoneskel}
The number of standard monomials of $M_{n}^{(1)}$ (and hence the dimension of the  $\mathbb K$-vector space $S/M_n^{(1)}$) is given by
\[|\stan_n^1| = (2n-1)(n-1)^{n-1}.\]
\end{prop}

\begin{proof}
Suppose $(a_1, a_2, \dots, a_n)$ is the exponent vector of a standard monomial of $M_n^{(1)}$.  Then by definition each entry $a_i$ is strictly less than $n$, and at most one entry equals $n-1$.  

If no entry equals $n-1$ then every entry satisfies $0 \leq a_i \leq n-2$ so we have $(n-1)^n$ possibilities. If exactly one entry (say $a_i$) equals $n-1$ then every other entry satisfies $0 \leq a_j \leq n-2$ for $j \neq i$.  Hence we have $n(n-1)^{n-1}$ possibilities.  Adding these up we get
\[(n-1)^n + n(n-1)^{n-1} = (2n-1)(n-1)^{n-1}\]
standard monomials, as desired.
\end{proof}

It turns out the number of elements in $\stan_n^1$ has a determinantal interpretation analogous to the case of classical parking functions.  Recall that $|\stan_n^{n-1}|$ is given by $\det \tilde {\mathcal L}(K_{n+1})$, the determinant of the reduced Laplacian of $K_{n+1}$.  For the one-skeleton a different but related matrix makes an appearance.

\begin{defn}
For a graph $G$ on vertex set $\{0,1, \dots, n\}$ the \defi{signless Laplacian} ${\mathcal Q} = {\mathcal Q}(G)$ is the symmetric $(n+1) \times (n+1)$ matrix with rows and columns indexed by the vertices of $G$ and with entries given by 
\[{\mathcal Q}_{i,j}= \begin{cases} \deg(i) &\mbox{if } i=j, \\ 
|\{\text{edges connecting $i$ and $j$}\}| & \mbox{if } i \neq j. \end{cases}\]
Define the \defi{reduced signless Laplacian} $\tilde {\mathcal Q}$ to be the matrix obtained from ${\mathcal Q}$ by deleting the row and column corresponding to the vertex 0.
\end{defn}

Note that ${\mathcal Q}$ has entries given by the absolute values of the entries of the usual Laplacian ${\mathcal L}$ (hence the name). For example, if $G=K_4$ we get the following matrices.

\[{\mathcal Q} = \begin{bmatrix}
3 & 1 & 1 &1 \\
1 & 3 & 1 & 1 \\
1 & 1 & 3 & 1\\
1 & 1 & 1 & 3
\end{bmatrix}
\hspace{.5 in}
\tilde {\mathcal Q }= \begin{bmatrix}
3 & 1 & 1 \\
1 & 3 &  1 \\
1 & 1 & 3 
\end{bmatrix}\]
\noindent
In this case one has $\det \tilde {\mathcal Q} = 20 = (5)(2^2)$ and in fact more generally we have the following.

\begin{prop}\label{prop:signless}
The number of standard monomials of $M_{n}^{(1)}$ is given by 
\[|\stan_n^1| = \det \tilde Q(K_{n+1}). \]
\end{prop}

\begin{proof}
According to Proposition \ref{prop:standardoneskel} it is enough to show that
\[\det \tilde{\mathcal Q}(K_{n+1}) = (2n-1)(n-1)^{n-1}.\]
\noindent
For this we examine the eigenvalues of the matrix $\tilde Q = \tilde Q(K_{n+1})$, which are also given (without proof) in OEIS $\#A176043$ \cite{OEIS}.  We have one eigenvalue $2n - 1$ with multiplicity 1 corresponding to the all 1's vector ${\bf 1}$.  Subtracting the matrix $(n-1)I_n$ from $\tilde Q$ gives us the matrix $J$ consisting of all ones, which has an $(n-1)$-dimensional kernel.  Hence $\tilde Q$ has one other eigenvalue $n-1$ with multiplicity $n-1$.  The result follows.
\end{proof}

Proposition \ref{prop:signless} suggests there may be a bijective proof of Proposition \ref{prop:standardoneskel} that extends well-known bijections between spanning trees and parking functions (see for instance \cite{ChePyl}).  For this we recall the following graph-theoretical interpretation of $\det \tilde {\mathcal Q}$.

\begin{prop}[\cite{Bap}, Theorem 7.8] \label{prop:TU}
For any graph $G$ the determinant of $\tilde {\mathcal Q}_G = \tilde {\mathcal Q}$ is given by 
\[\det \tilde {\mathcal Q} = \sum_H 4^{c(H)}, \]
where the summation runs over all spanning $TU$-subgraphs $H$ of $G$ with $c(H)$ unicyclic components, and one tree component which contains the vertex 0.
\end{prop}

Here a \defi{unicylic} graph is a graph with a single cycle.  A \defi{$TU$-subgraph} is a subgraph of $G$ whose components are trees or unicylic graphs with odd cycles.  

\begin{rem}
The signless Laplacian also has connections to the theory of \emph{signed graphs} as developed by Zaslavsky in \cite{Zas}. For any graph $G$ the signless Laplacian ${\mathcal Q}(G)$ can be factored as ${\mathcal Q}(G) = j \circ j^T$, where $j$ is the \emph{signless incidence matrix} of $G$.  The matrix $j$ can be taken as an incidence matrix for the underlying \emph{signed graph} of $G$, where all edges are taken to be negative. From \cite{Zas} it is known that the independent sets in the linear matroid determined  by the matrix $j$ correspond to sets of edges where each component either contains no circles, or just one circle which is `negative' (in this context equivalent to having an odd number of vertices).  Applying Cauchy-Binet to calculate $\det \tilde{\mathcal Q}$ involves computing determinants of maximal minors of $\tilde j$, and the nonzero contributions are exactly those subgraphs of $G$ described above.
\end{rem}

There are a number of combinatorial bijections between the set ${\mathcal P}(n)$ (which we identify with $ \stan_n^{n-1}$) and the set of spanning trees of the complete graph $K_{n+1}$, including Dhar's burning algorithm and other variations (\cite{ChePyl}). Many of these extend to the context of arbitrary graphs. It is a natural question to find a similar map between the spanning $TU$-subgraphs of $K_{n+1}$ and the set $\stan_n^{n-2}$.  We note that a spanning tree is, in particular, a spanning $TU$ subgraph and hence the desired map should extend one of these bijections.

A bijective proof of Proposition \ref{prop:standardoneskel} would associate to each spanning $TU$-subgraph $H \subset K_{n+1}$ a collection of $4^{c(H)}$ elements of $\stan_n^1$.  Each spanning tree of $K_{n+1}$ would be assigned $4^0 = 1$ standard monomials, so presumably such a bijection would extend the correspondence between usual parking functions and spanning trees.  In Example \ref{ex:K4oneskel} we have the 16 parking functions coming from the spanning trees of $G$, and 4 new standard monomials coming from the $TU$-subgraph consisting of the edges $\{12,13,23\}$.  It would also be interesting to see if the degree sequence of the elements in $\stan_n^1$ can be related to a notion of `inversion' for the set of spanning $TU$-subgraphs (see Section \ref{sec:TUinv}).

\section{Arbitrary graphs and further questions} \label{sec:conclusion}

As we have seen, the study of $\stan_n^k$ for various values of $k$ leads to combinatorial notions that relate and extend some known interpretations of classical parking functions. As mentioned in the introduction, the notion of a parking function can be generalized to $G$-parking functions, where $G$ is any graph on vertex set $V = \{0,1,\dots, n\}$ (with specified root vertex $0$).  Here we fix such a graph $G$ and for a subset $\sigma \subset [n]$ we let $\deg_\sigma(i) = |\{j \in V(G) \backslash \sigma: i \sim j \}|$ denote the number of vertices adjacent to $i$ that are outside $\sigma$.  For any subset $\sigma \subset [n]$ define a monomial
\[m_\sigma = \prod_{i \in \sigma}x_i^{\deg_\sigma(i)}.\]
By definition the \defi{$G$-parking functions} correspond to the standard monomials of the ideal $M_G \subset S$ generated by all such $m_\sigma$ for $\sigma \neq \emptyset$. We define the ideal $M_G^{(k)}$ to be the subideal of $M_G$  generated by monomials corresponding to subsets of size at most $k+1$.   A natural question to ask is if the objects studied here can be extended to the context of general graphs.  

To motivate this study, note that for an arbitrary graph $G$ (with specified root vertex), the standard monomials of the $k$-skeleton ideal $M_G^{(k)}$ specialize to two natural generalizations of classical parking functions: a certain class of ${\bf u}$-parking functions on the one hand (by taking $G=K_{n+1}$), and $G$-parking functions on the other (by taking $k = n-1$).   In this context one could hope for a generalization of Corollary \ref{cor:standard} that counts standard monomials of $M_G^{(k)}$, this time incorporating data from the underlying graph $G$.

\subsection{Codimension-one ideals.}
As we have seen, the elements of $\partial(\stan_n^{n-2}) = \stan_n^{n-2} \backslash \stan_n^{n-1}$ naturally correspond to the set of \emph{prime} parking functions.  It is not clear if there is a similar interpretation for arbitrary $G$.  For example it does not seem that elements of $\partial(\stan_G^{n-2})$ correspond in any natural way to a subset of usual $G$-parking functions. 

Computing the number of elements in $\partial(\stan_G^{(n-2)})$ for arbitrary $G$ seems like a difficult task.  In \cite{KumLatSon} the authors determine the size of these sets for the case of $G = K_{n+1} \backslash \{e\}$ for some edge $e$. In the case that $e$ is adjacent to the root vertex $0$ they show that 
\[|\partial(\stan_G^{(n-2)})| = (n-1)^{n-1},\] whereas if $e$ is not adjacent to the root we have
\begin{equation}\label{eq:completeminusedge}
    |\partial(\stan_G^{(n-2)})| = (n-1)^{n-3}(n-2)^2.
\end{equation}

We have seen that the elements of $\partial(\stan_n^{(n-2)})$ are counted by uprooted trees on $[n]$ and it is not clear if $\partial(\stan_G^{(n-2)}))$ can be related to some subset of the spanning trees of $G$.  We remark that the proof of Equation \ref{eq:completeminusedge} from \cite{KumLatSon} involves showing that the set  $\partial(\stan_G^{(n-2)})$ is in bijection with the set of uprooted trees on $[n]$, where $1$ is not adjacent $n$. 

Also, Gessel \cite{Ges} has generalized the notion of inversion in the context of \emph{$\kappa$-inversions} for spanning trees of an arbitrary graph $G$. In \cite{PerYanYu} Perkinson, Yang, and Yu give a bijection between $G$-parking functions and spanning trees that preserves degree and the number of $\kappa$-inversions.  Hence one could search for a generalization of Theorem \ref{thm:inversions} that extends to $\partial(\stan_G^{n-2})$.

\subsection{One-dimensional ideals.}

The results from Section \ref{sec:signless} can also be considered in the context of arbitrary graphs. Suppose $G$ has vertex set $\{0,1, \dots, n\}$ and as above let $\tilde {\mathcal Q}_G$ denote its reduced signless Laplacian. Our convention here is that $0$ is taken to be the sink (corresponding to which row/column should be deleted) but note that in general $\det(\tilde {\mathcal Q}_G)$ depends on this choice.   This is in contrast to the usual Laplacian, where the determinant simply counts the number of spanning trees containing the sink (and hence is independent of choice of sink).   In fact the natural extension of Proposition \ref{prop:signless} does not hold for general $G$, as the following example illustrates. 

\begin{example}\label{ex:K5minusedge}
Let $H$ be the graph obtained from removing the edge $(34)$ from the graph $K_5$.  The reduced signless Laplacian is 
\[\tilde {\mathcal Q}_H = \begin{bmatrix}
4 & 1 & 1 & 1 \\
1 & 4 & 1 & 1 \\
1 & 1 & 3 & 0\\
1 & 1 & 0 & 3
\end{bmatrix}
\]
\noindent
with $\det(\tilde {\mathcal Q}_H) = 99$.  According to Macaulay2 \cite{M2} we have $\dim_{\mathbb K} M_G^{(1)} = 105$, so that there are $105$ standard monomials of $M_G^{(1)}$ in this case. Note however if $H^\prime$ is the graph obtained from $K_5$ by removing the edge $01$ then we get $\det(\tilde {\mathcal Q}_{H^\prime}) = 135$, while there are 135 standard monomials of $M_{H^\prime}^{(1)}$.
\end{example}

\begin{figure} \label{fig:K5minusedge}
\includegraphics[scale = .4]{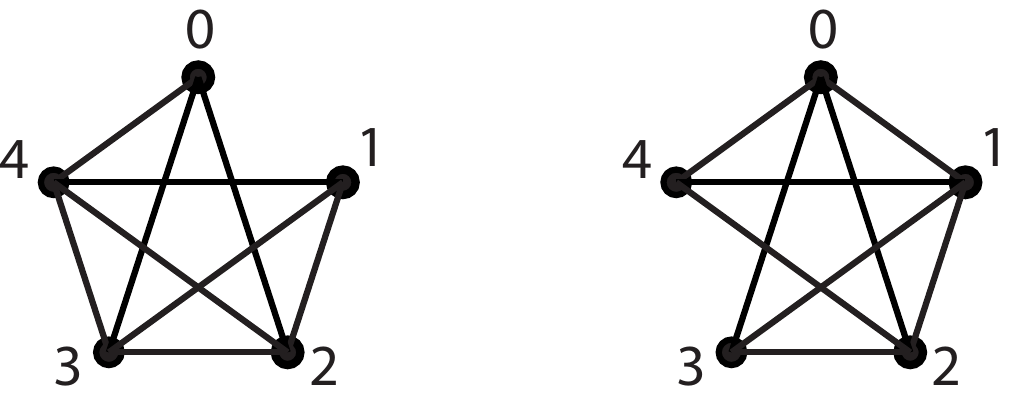}
\caption{The graphs $H^\prime$ and $H$ from Example \ref{ex:K5minusedge}.  The graphs are isomorphic but have different values of $\det \tilde {\mathcal Q}$ (in both cases the sink is given by the vertex 0).  }  
\end{figure}

After a number of calculations we have not found an example where the determinant of $\tilde {\mathcal Q}$ is \emph{larger} than the dimension of $S/M_G^{(1)}$, which begs the following question.

\begin{question} \label{ques:inequality}
For any graph $G$ is it true that
\[|\stan(M_G^{(1)})| \geq \det(\tilde {\mathcal Q}_G)?\]

\end{question}
\noindent

\begin{rem}
After a version of this paper was posted to the arXiv this question was answered in the affirmative by Kumar, Lather, and Roy \cite{KumLatRoy} in the more general setting of multigraphs.
\end{rem}

\subsection{Inversions in TU-subgraphs}\label{sec:TUinv}

As we have seen, the degree sequences of elements in $\stan_n^k$ can often be interpreted as statistics on related combinatorial objects. For any $k$ we define a generating function 
\[P^{(k)}_n(q) = \sum_{\beta = (b_1, b_2, \dots, b_n)} q^{b_1 + b_2 + \dots +b_n}\]
where $\beta$ ranges over all standard monomials of $M_n^{(k)}$.

For the case of $k=n-1$ and $k=n-2$ the coefficients of these polynomials can be related to the number of inversions in associated labeled tree-like structures.

In Section \ref{sec:signless} we have seen that elements of $\stan_n^1$ are counted (in a weighted fashion) by the spanning $TU$-subgraphs of $K_{n+1}$ (with spanning trees appearing as a subclass).  A natural question to ask is whether the degree of elements in $\stan_n^1$ can be interpreted by `inversions' in a similar way. For small values of $n$ we have
\begin{equation}
\begin{split}
& P^{(1)}_1(q) = 1 \\
& P^{(1)}_2(q) = 2q + 1\\
& P^{(1)}_3(q) = 3q^4 + 7q^3 + 6q^2 + 3q + 1\\
& P^{(1)}_4(q) = 4q^9 + 13q^8 + 28q^7 + 38q^6 + 40q^5 + 31 q^4 + 20 q^3 + 10 q^2 + 4q +1
\end{split}
\end{equation}

\subsection{Other skeleta and trees}

Another open question here is whether or not standard monomials of other skeleton ideals are related to other tree-like structures. Classically we have that $\stan_{n}^{n-1}$ is in bijection with the number of labeled trees on $[n]$ (equivalently labeled rooted forests on $[n-1]$).  From our work we see that $\stan_n^{n-2}$ counts labeled trees and uprooted labeled trees, whereas $\stan_n^{1}$ counts labeled trees and spanning $TU$-subgraphs (a generalization of spanning trees related to basis of an underlying signed graph).  We do not know if $\stan_n^{k}$ is related to other tree-like structures for other values of $k$.

\subsection{Other subcomplex ideals}
In this paper we studied skeleta of $G$-parking function ideals, but the constructions makes sense in a more general setting.  For instance if $\Delta$ is any simplicial complex on vertex set $[n]$ we can study the ideal $M_n^\Delta$ defined in the natural way
\[M^\Delta_n = \langle m_\sigma : \sigma \in \Delta\rangle,\]

\noindent
where the $m_\sigma$ for nonempty $\sigma \subset [n]$ are the monomials defined above.

Backman considered this level of generality in his \emph{hereditary chip-firing} models in \cite{Bac}.  One can then ask for a formula for  $|\stan_n^\Delta|$, the number of standard monomials of this ideal, perhaps in terms of some statistic on $\Delta$.  For instance does the number of standard monomials only depend only on certain combinatorial features of $\Delta$? Does the topology play any role? Perhaps one can obtain a nice formula for $|\stan_n^\Delta|$ when $\Delta$ is a simplicial complex with desirable combinatorial properties (for example matroidal, shellable, shifted, etc.).

{\bf Acknowledgements.}  We were introduced to the ideals $M_G^{(k)}$ by Spencer Backman, who suggested the idea of studying algebraic and combinatorial properties of subideals of $M_G$ generated by restricted subsets.  We thank Spencer and also Robin Sulzgruber, Suho Oh, and Catherine Yan for valuable discussions. We are grateful to the anonymous referees who provided helpful comments on this and earlier versions of the paper.

\end{document}